\documentclass[10pt,journal,compsoc]{IEEEtran}

\ifCLASSOPTIONcompsoc
  \usepackage[nocompress]{cite}
\else
  \usepackage{cite}
\fi
\ifCLASSINFOpdf
   \usepackage[pdftex]{graphicx}
\else
   \usepackage[dvips]{graphicx}
\fi
\usepackage{amsmath}
\usepackage{nccmath}
\usepackage{algorithm}
\usepackage{algorithmicx}
\usepackage{algpseudocode}

\usepackage{array}

\ifCLASSOPTIONcompsoc
 \usepackage[caption=false,font=footnotesize,labelfont=sf,textfont=sf]{subfig}
\else
 \usepackage[caption=false,font=footnotesize]{subfig}
\fi
\usepackage{xcolor}
\definecolor{edit_col}{rgb}{0, 0, 0}

\begin{document}
%
\title{Backpropagation on Dynamical Networks}
%
%
%
%

\author{Eugene Tan,
        D\`{e}bora Corr\^{e}a,
        Thomas Stemler,
        Michael Small
\IEEEcompsocitemizethanks{\IEEEcompsocthanksitem E. Tan, T. Stemler and M. Small are with the Department of Mathematics \& Statistics, The University of Western Australia, Crawley, WA 6009, Australia\protect\\
E-mail: eugene.tan@uwa.edu.au
\IEEEcompsocthanksitem D. Corr\^{e}a is with the Department of Computer Science \& Software Engineering, The University of Western Australia, Crawley, WA 6009.}
\thanks{Manuscript received INSERT DATE; revised INSERT DATE.}}

%
%

\markboth{IEEE TRANSACTIONS OF NETWORK SCIENCE AND ENGINEERING}%
{IEEE TRANSACTIONS OF NETWORK SCIENCE AND ENGINEERING}
%



\IEEEtitleabstractindextext{%
\begin{abstract}
 Dynamical networks are versatile models that describe a variety of behaviours such as synchronisation and feedback. However, applying these models in real systems is difficult as prior information of the connectivity structure or local dynamics is often unknown and must be inferred from observations of network states. Additionally, the influence of coupling interactions further complicates the isolation of local node dynamics. Given the architectural similarities between dynamical networks and recurrent neural network (RNNs), we propose a network inference method based on the backpropagation through time algorithm used to train RNNs. This method aims to simultaneously infer both the connectivity structure and local node dynamics from node state observations. An approximation of local node dynamics is first constructed using a neural network. This is alternated with an adapted BPTT algorithm to regress corresponding network weights by minimising prediction errors of the network based on the previously constructed local models until convergence. This method was found to be successful in identifying the connectivity structure for coupled networks of chaotic oscillators. Freerun prediction performance with the resulting local models and weights was comparable to the true system with noisy initial conditions. The method is also extended to asymmetric negative coupling.
\end{abstract}

}

\maketitle

\IEEEdisplaynontitleabstractindextext

%
\IEEEpeerreviewmaketitle

\IEEEraisesectionheading{\section{Introduction}\label{sec:introduction}}
Dynamical networks are a common occurrence when studying complex systems where the goal is to model large systems of multiple interacting components. In its simplest form, a dynamical network may be defined as having three main components: local node dynamics, coupling dynamics and connectivity structure (see Equation \ref{eq:DynamicalNetwork}). One key feature of dynamical networks is its ability to recreate and describe rich and interesting dynamics such as chimera states, synchronisation \cite{abrams2004chimera,andreev2019chimera, zhu2014chimera} and cascading failure commonly encountered in real world systems \cite{majhi2019chimera}. They have also been applied to other systems such as neuron networks \cite{izhikevich2007dynamical, majhi2019chimera}, power grids \cite{wu2010evolution, liu2010exponential}, epidemic spread \cite{hota2021impacts}, and cardiac arrhythmia \cite{rosenblum2019dynamical}.

\begin{equation}
\dot{x}_i = f(x_i) + \sum_{i\neq j}c_{i,j}g(x_i,x_j).
\label{eq:DynamicalNetwork}
\end{equation}

In Eq \ref{eq:DynamicalNetwork}, nodes are indexed with $i$, $c$ represents the coupling structure, $f$ and $g$ represents the local and coupling dynamics respectively. Despite its relatively simple form, data-driven applications of this framework to model real world systems can be difficult if information regarding the coupling structure and dynamics is not accessible. In these cases, these components need to be inferred before a dynamical networks framework can be applied. One simple example would be the data driven analysis of a network of chaotic oscillators. Specifically, there is the inverse problem of inferring network characteristics from data. However, this problem is difficult due to the complicated interplay between local models and coupling behaviour. \textcolor{edit_col}{Many existing methods of inferring dynamical network properties require either the local or coupling dynamics to be at least partially known. Alternatively, some methods restrict inferences to statistical arguments with uncertainty in order to remove the dependence on a priori information.}

Numerous approximation and inference methods for inferring network connectivity have been developed each with varying success. These methods may be classified into three main approaches:

\begin{enumerate}
    \item \textbf{Direct Approach} - For systems where both local $f$ and coupling dynamics $g$ is explicitly known, a direct algebraic approach using derivatives have been used to recover connectivity structure \cite{shandilya2011inferring}.
    \item \textbf{Perturbation Methods} - These methods try and recreate a proxy dynamical network that performs similarly to the target system. Small perturbations are applied to individual nodes in the network and observed to see how each signal propagates. These propagations are then used to infer the network connectivity structure \cite{stepaniants2020inferring, banerjee2019using, banerjee2021machine} .
    \item \textbf{Signal Correlation} - This approach focuses on quantifying the correlation between node signals \cite{napoletani2008reconstructing}. Node pairs with highly correlated signals are used to infer the connectivity structure. Maximum route entropy has also been used as an alternative to achieve the same result. \cite{weistuch2020inferring}. An extension of this approach is the employment of some notion of causality such as Granger causality \cite{paluvs2019coupling, kornilov2020reconstruction} and short term causal dependence. \cite{stepaniants2020inferring}. 
\end{enumerate}

The inference of connectivity structure also poses an additional challenge of false positive couplings and remains an open research area. Granger causality and similar statistical based approaches have been found to particularly struggle with this problem, even for simple network structures such as rings and chains \cite{kornilov2020reconstruction, wu2012inferring}. 

The second problem of constructing models to describe the local dynamics $f$ in multi-node networks has been studied for pseudo-periodic systems \cite{rosenblum2019dynamical,paluvs2019coupling}. Outside these, many approaches focus on the case of a single node without external interactions and utilise a time-delay embedding with appropriately chosen delay and embedding dimension \cite{takens1981detecting}. Proposed modelling methods for long term prediction include neural network predictors \cite{zhang2000predicting}, radial basis predictors \cite{small2002minimum} and reservoir computing \cite{jaeger2001echo} among several others. Assessing the quality of these local models \textcolor{edit_col}{relies} on the comparison of dynamical invariants such as Lyapunov exponents or correlation dimension \cite{pathak2017using,lu2018attractor,haluszczynski2019good} and attractor homology \cite{tan2021grading}. 

\textcolor{edit_col}{There is also a third related problem of identifying coupling equations $g$ from signals. Solutions have been proposed using phase dynamics for networks exhibiting pseudo-periodic behaviour, but not necessarily for aperiodic behaviour \cite{rosenblum2019dynamical}. For simple oscillator networks, Panaggio et al.  \cite{panaggio2019model} propose a gradient descent method to regress coupling functions of a given form. This approach is draws similarities with \cite{shandilya2011inferring}. For simplicity, this paper assumes that the coupling equations are known.}

\textcolor{edit_col}{For low dimensional systems, symbolic approaches aimed at tackling either the dynamics or connectivity structure inference problem has gained increased traction in recent times. A method based on the ideas of compressed sensing \cite{wang2016data} and sparse optimisation methods (SINDy) \cite{brunton2016discovering} have been proposed to reconstruct the functional form of the dynamics equations using a library of basis functions (e.g. power series, Fourier series) for isolated dynamical systems with no interactions. Compressive sensing tries to find the sparsest combination of coefficients for some power series basis via convex optimisation with respect to the $L_1$-norm. Link and coupling weight identification can also be included into this approach \cite{wang2016data, wang2011time}. However, we note that this method may be limited to simple systems with analytical or well behaved equations.}

\textcolor{edit_col}{The symbolic approach with basis functions was also extended to local dynamics and coupling functions by utilising prior known information of the network connectivity structure \cite{gao2022autonomous}. A similar approach (ARNI) \cite{casadiego2017model} using basis functions has also been proposed to reconstruct the coupling function and connectivity structure of dynamical networks. ARNI is also applicable to hetereogeneous networks where local dynamics $f$ differ between nodes. We note that ARNI can also be adapted to infer local dynamics instead \cite{gao2022autonomous}.}

\textcolor{edit_col}{Despite these topics being well studied, the literature on inferring both connectivity structure and node dynamics from time series observations is not rich. Methods for reconstructing local dynamics (e.g. SINDy \cite{brunton2016discovering} and Gao \cite{gao2022autonomous}) are either restricted to single node systems, or rely or prior information about the connectivity structure. In cases where the aim is to infer connectivity structure, either local dynamics are known prior \cite{shandilya2011inferring} or are not directly inferred \cite{casadiego2017model}. 
There remains the need of a method that tackles the inference of both of these properties simultaneously.}

\textcolor{edit_col}{We note that a method utilising a mean field approach has been proposed to tackle the simultaneous inference of local dynamics and connectivity structure \cite{eroglu2020revealing}. This method constructs an `effective network' that acts as an approximate proxy for the real dynamical network, reproducing many of its properties and behaviours. Under specific assumptions such as weak coupling, network heterogeneity and limited links, this method has been shown to perform well in identifying critical transitions and the connectivity structure of a cat cortex. However, the `effective network' relies on the network sparsity and presence of hubs to recover an approximation of the local dynamics. Additionally, this approximation of local dynamics is not fully disentangled from the interaction dynamics.}

In this paper, we present a proof of principle of a data-driven approach for simultaneously inferring both the local node dynamics and weighted connectivity structure from the observation of node signals, assuming only prior knowledge of the coupling function. Recognising the structural similarities between dynamical networks and RNNs, we propose an adaptation of the backpropagation through time (BPTT) algorithm used to train RNNs \cite{werbos1990backpropagation}, as a method of regressing the connectivity structure of a dynamical network. This regression is then used to decouple the observed node signals and simultaneously construct a model of local dynamics. \textcolor{edit_col}{This method combines the tasks achieved by several of the abovementioned methods, whilst also dealing instances of dynamical networks with weighted edges.}

We test the performance of the proposed method on a dynamical network of Lorenz \cite{lorenz1963deterministic} and Chua oscillators \cite{kengne2017dynamics} with a randomly initialised weighted connectivity structure. We find that an adapted backpropagation regression method was able to reproduce the network coupling structure and local dynamics for a 64 node network. Additional cases with asymmetric and negatively weighted coupling are also explored, yielding similar results. To test the performance on a more realistic application, the proposed backpropagation regression method was applied to a network of driven FitzHugh-Nagumo oscillators operating in the chaotic regime to simulate a simple biological neuron network. Similar to the previous cases, the backpropagation method was able to successfully reproduce both the network structure and local dynamics. 

\section{Background}

\subsection{Dynamical Networks}

Dynamical networks can be defined as a network $G=(f,g,C)$ consisting of three main components, local dynamics $f$, coupling dynamics  $g$ and connectivity structure $C$. Node states evolve according to Equation \ref{eq:DynamicalNetwork} with local behaviour governed by $f$ and additional external influences from neighbours due to coupling components $g$ and $C$.

The effect of coupling between nodes results in a high dimensional system with multiple interacting sub-systems. For small or highly connected networks, sufficiently large coupling can cause network-wide synchronisation of node states. As coupling is decreased, chimera states and group synchronisation emerges \cite{abrams2004chimera,andreev2019chimera, zhu2014chimera}. When no coupling is present, all nodes behave independent of each other.

Here, we focus on dynamical networks based on two different 3-dimensional chaotic systems, Lorenz and Chua, with diffusive coupling in the first component. For the latter, we will focus on the Chua system with cubic nonlinearity $\phi$ operating in the single scroll mode \cite{kengne2017dynamics} to analyse dynamics contrasting with the Lorenz system.

\subsection{Backpropagation and Recurrent Neural Networks}
Recurrent neural networks (RNNs) are a modelling architecture developed in the field of machine learning that is widely used for sequence and time series prediction \cite{raubitzek2021taming, dubois2020data, Sangiorgio2020}. RNNs consist of an input layer with weights $C^{(in)}$, a fully connected hidden layer of nodes $\mathbf{h}(t)$ with feedback weights $C^{(h)}$, and an output layer with readout weights $C^{(out)}$. For prediction, an input sequence $\mathbf{x}(t_n)$ is fed into the RNN via the input layer with \textcolor{edit_col}{input weights $C^{(in)}$. These inputs are usually discretely sampled at regular time intervals ${t_n}$ and then sequentially fed into the recurrent hidden layer whose nodes feedback to each other according to weights before applying a nonlinear transformation $\sigma$. In its continuous time formulation, this allows network to approximate nonlinear dependencies between states with nodes evolving according to Equation \ref{:RNN},}

\begin{equation}
    \mathbf{\hat{h}}(t) = \sigma (C^{(in)}\mathbf{x}(t) + C^{(h)}\mathbf{h}(t)+\mathbf{b}).
    \label{:RNN}
\end{equation}

Predictions $\hat{\mathbf{x}}(t)$ are calculated using the readout weights and hidden node states as given in Equation \ref{:Readout},
\begin{equation}
    \mathbf{\hat{x}}(t) = C^{(out)}\cdot\mathbf{h}(t).
    \label{:Readout}
\end{equation}
The training of RNNs involve the small adjustment of network weights $C^{(in)}$, $C^{(h)}$, $C^{(out)}$ and biases $\mathbf{b}$. This is achieved via the backpropagation through time (BPTT) algorithm \cite{paluvs2019coupling}. The BPTT algorithm can be interpreted as an adaptation of the gradient descent training algorithm for neural networks that accounts for the propagation of network weight effects through the history of the RNN's hidden node states.

\begin{figure}
    \centering
    \includegraphics[width = 0.35\textwidth]{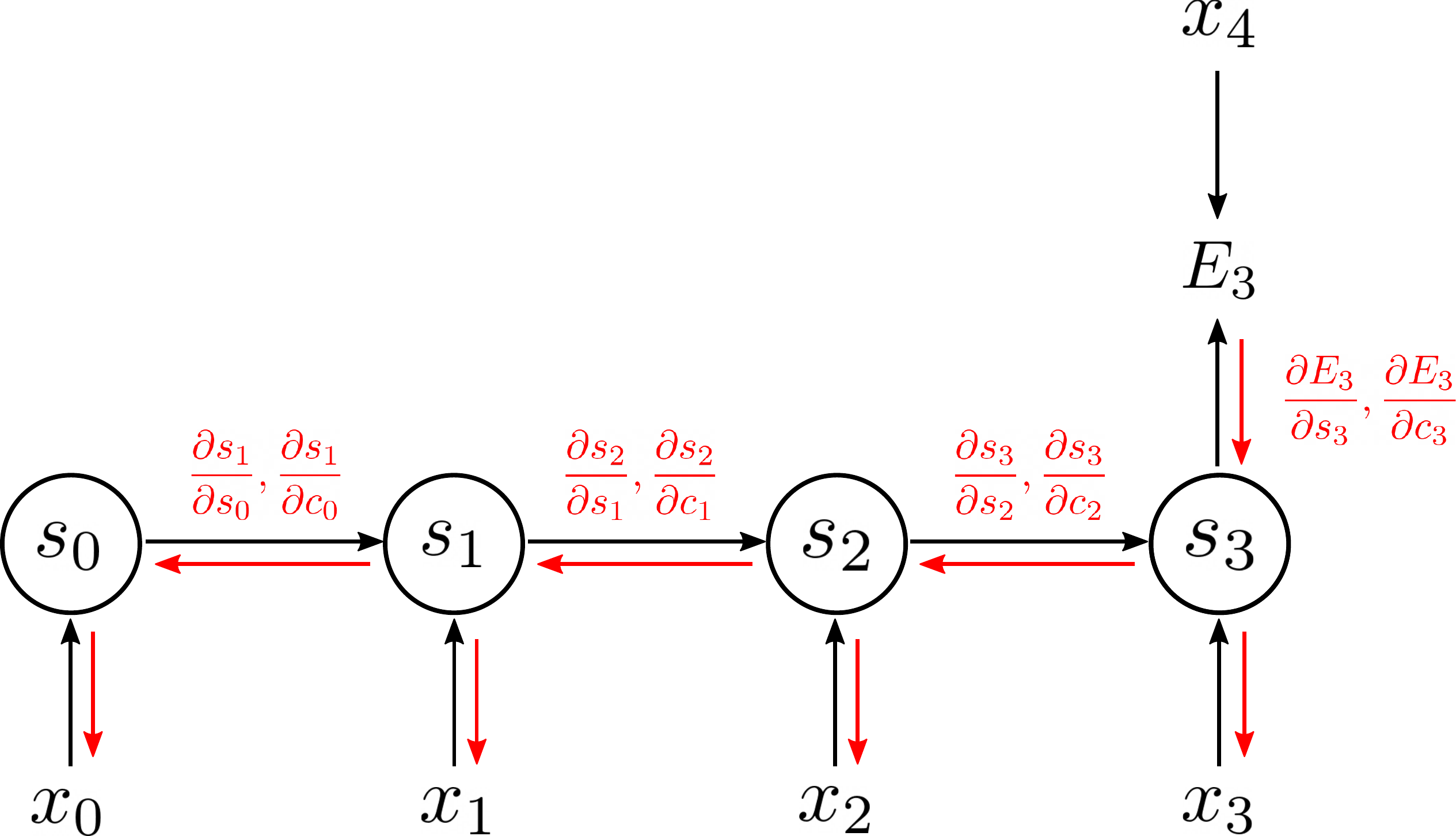}
    \caption{Schematic of typical backpropagation through time calculation for training recurrent neural networks with an input length of 4. Input time series $x_t$ is used to evolve node states $s_t$ during the forward pass. Node states are used to calculate a value for the error $E_t$ at a given time $t$. \textcolor{edit_col}{Derivatives are taken with respect to node states and weights. Variables $c_n$ correspond to the weight of the link joining $s_{n-1}$ to $s_n$ and is contained within $C^{h}(t)$}}.
    \label{fig:NormalBackprop}
\end{figure}

\textcolor{edit_col}{The BPTT algorithm works by unfolding the predictions of the recurrent network at each time step. This process reframes the forward pass of an RNN into a feedforward network where each layer corresponds to the node states of the RNN at a specific time.} For a given RNN model, an input sequence $x_{in}=\{ x_{0},x_{1},...,x_{\tau} \}$ is fed in to produce a sequence of outputs $x'_{out}=\{ x'_{\tau+1},x'_{\tau+2},...,x'_{\tau+k} \}$. Prediction losses are calculated with respect to some metric $d$ and the real output $x_{out}=\{ x_{\tau+1},x_{\tau+2},...,x_{\tau+k} \}$,
\begin{equation}
    \mathcal{L}=d(x'_{out}, x_{out}).
\end{equation}
The gradient of the loss taken with respect to node weights is calculated by propagating the error backward through time steps and used to update node weights (see. Figure \ref{fig:NormalBackprop}).

From observation, the hidden layer of an RNN functions identically to a dynamical network. We aim to investigate if the BPTT approach can be adapted to regress the weights of a dynamical network given a set of input node states.

\section{Methods}

\begin{figure}
    \centering
    \includegraphics[width = 0.3\textwidth]{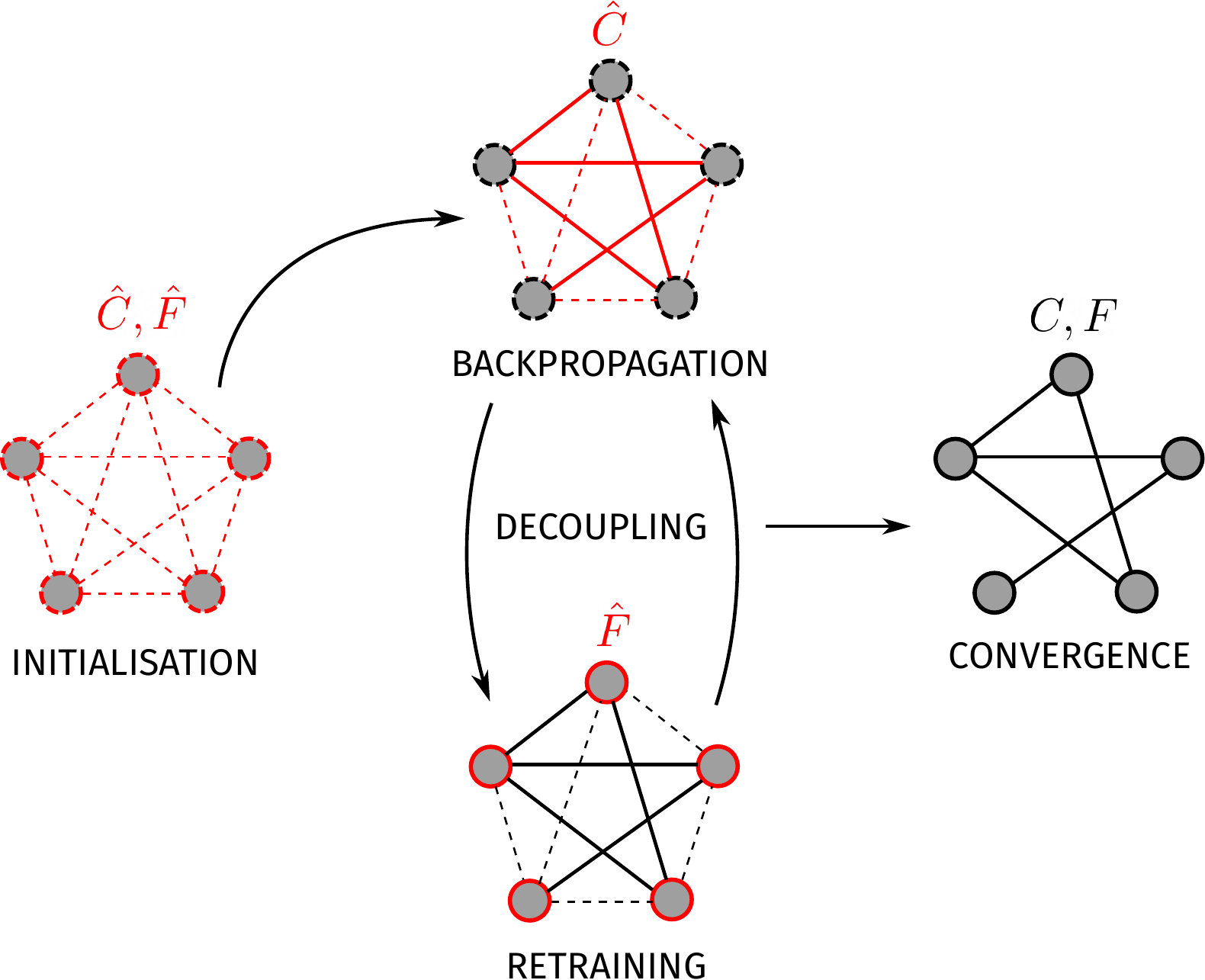}
    \caption{Overview of the steps in the proposed backpropagation regression method. Components being trained (red) in each step. Regressed components are given with solid lines.}
    \label{fig:Flowchart}
\end{figure}

This paper proposes a novel method of simultaneously inferring the local dynamics $f$ and connectivity structure $C$ of a dynamical network with a known coupling function $g$ by only observing node states. The proposed method consists of three main stages, namely initialisation, backpropagation and decoupling. A flowchart outlining this process is provided (see Figure \ref{fig:Flowchart}). 


\subsection{Backpropagation Regression Algorithm}

\subsubsection*{Initialisation}

The initialisation stage consists of constructing a local model $\hat{f}$ that approximates the true local dynamics $f$. Due to the complex interplay between the local $f$ and coupling $g$ dynamics, getting an accurate coupling-free estimation of $f$ directly from node signals is not possible. 

\begin{figure}
    \centering
    \includegraphics[scale = 0.15]{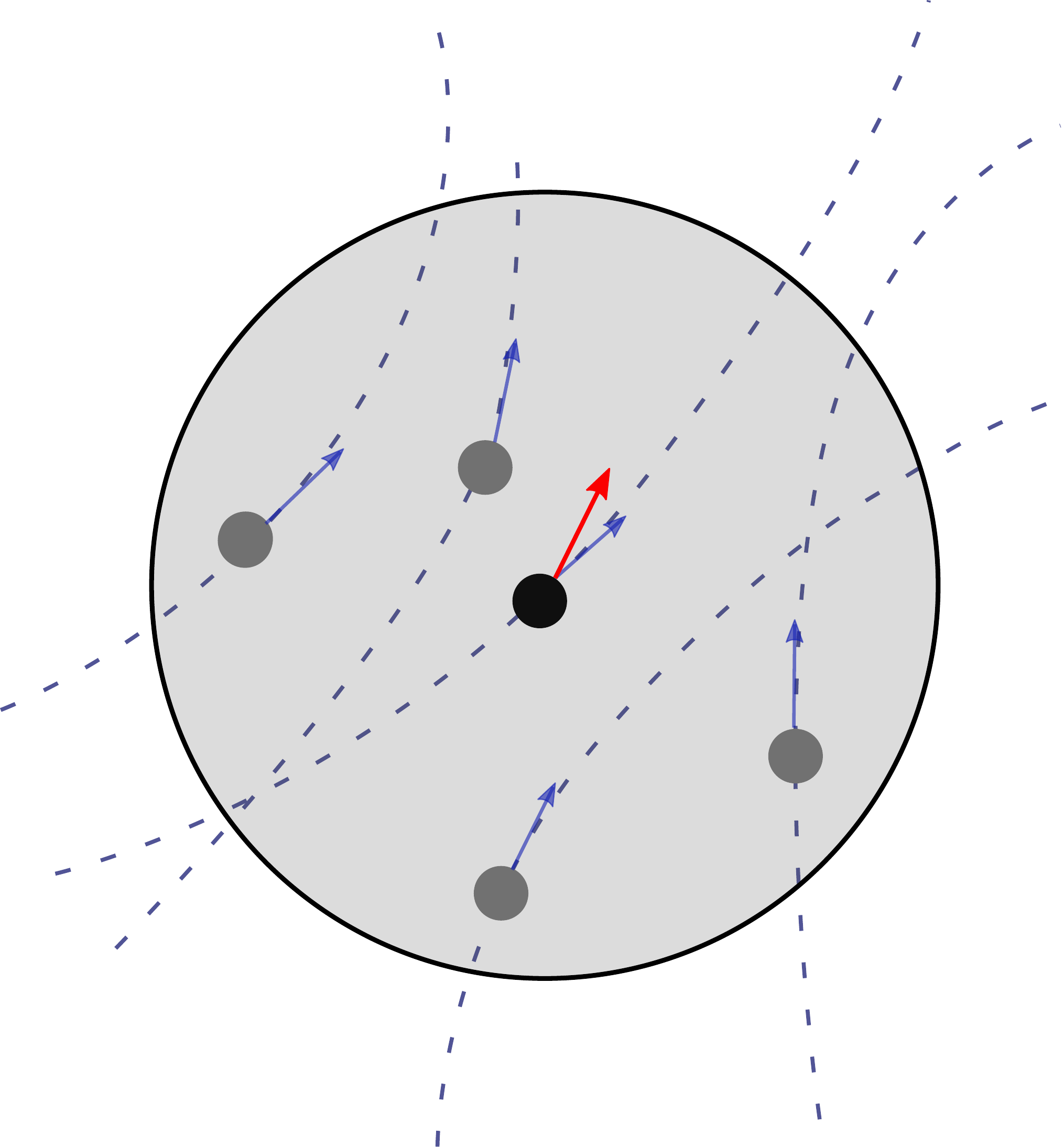}
    \caption{Mean field approach where effective vector field (blue) at a point in phase space is taken as the average of \textcolor{edit_col}{$K$} nearest neighbours' vector field influenced with neighbouring node interactions.}
    \label{fig:MeanField}
\end{figure}

The estimated local model $\hat{f}$ is defined as a map from state space $\mathbf{x}(t)$ to vector field $\mathbf{\dot{x}}(t)$. The time series is numerically differentiated in order to retrieve $\mathbf{\dot{x}}(t)$,

\begin{equation}
    \mathbf{\dot{x}}(t) \approx \frac{\mathbf{\dot{x}}(t+\delta t)-\mathbf{\dot{x}}(t)}{\delta t}.
\end{equation}
For this stage of the method, we assume that the coupling effects from $g$ is relatively small compared to local dynamics $f$. This allows a mean field approach to be taken when constructing an approximate local model $\hat{f}$. \textcolor{edit_col}{We note that this method requires that the `attractor' (coupling with other neighbours in the dynamical network may cause spurious intersections of trajectories) is not too distorted from the attractor of the non-coupled case.} Here, we calculate an estimate of the true vector field at each point in state space as an average over its $K$ nearest neighbours (see Figure \ref{fig:MeanField}),
\begin{equation}
    \mathbf{\hat{\dot{x}}}(t) = \frac{1}{K}\sum_{i=1}^K \frac{\mathbf{\dot{x}}^{(i)}(t+\delta t)-\mathbf{\dot{x}}^{(i)}(t)}{\delta t},
\end{equation}
\textcolor{edit_col}{where $\mathbf{\dot{x}^{(i)}}(t)$ corresponds to the velocity of the $i^{th}$ closest neighbour. The identification of neighbours are done with respect to the Euclidean distance with all points in the time series (i.e at all times $t$) being considered potential candidates of nearest neighbours.}

The estimated local dynamics $\hat{f}$ is defined as a mapping $\hat{f}:\mathbf{x}(t)\to\mathbf{\hat{\dot{x}}}(t)$. A feedforward network with 2 hidden layers \textcolor{edit_col}{was} used to learn the mapping due to its simplicity and ease of implementation. Each hidden layer is composed of 128 neurons with a \textcolor{edit_col}{sigmoid activation function.} This estimated local dynamics $\hat{f}$ calculated from the mean field approach acts as an initial proxy for the true local dynamics, which will then be further refined to remove coupling effects on the vector field during the backpropagation and retraining stages of the algorithm. In addition to the local dynamics $\hat{f}$, values for the estimated coupling adjacency matrix $\hat{C}$ must also be initialised. However, there is no restriction on what values of $\hat{C}$ may be taken.

\subsubsection*{Backpropagation}
The main aim of the backpropagation stage is to provide an improved estimate for the coupling adjacency matrix $\hat{C}$. Similar to the BPTT algorithm, a forward pass consisting of freerun predictions is first calculated using the estimated $\hat{f}$ and  $\hat{C}$ with randomly selected initial values from the observed time series in order to properly sample across the whole attractor,

\begin{subequations}
    \begin{align}
        \hat{\dot{\mathbf{x}}}_i = \hat{f}(\mathbf{x}_i) + \sum_{i\neq j}\hat{c}_{i,j}g(\hat{\mathbf{x}}_i,\hat{\mathbf{x}}_j),
        \\
        \hat{{\mathbf{x}}}_i(t+\delta t) = \hat{{\mathbf{x}}}_i(t) + \delta t \; \hat{\dot{\mathbf{x}}}_i(t).
    \end{align}
    \label{eq:forward}
\end{subequations}
\textcolor{edit_col}{For simplicity, we assume coupling only in the first component of each time series.}

Freerun predictions for all nodes $\mathbf{\hat{x}_i}(t)$ are calculated by recursively evaluating Equations \ref{eq:forward}. The loss in freerun predictions due to errors in the local model $\hat{f}$ and network weights $\hat{C}$ can be calculated using some distance function. Here, $n$ is the length of the freerun prediction, $N$ is the number of nodes. This loss at each node, $E_i$, is calculated using the $L_2$-norm across the total freerun trajectory \{$\hat{\mathbf{x}}_i(t_0),...,\hat{\mathbf{x}}_i(t_n)$\},

\begin{equation}
    \mathcal{L} = \sum_{t=t_0}^{t=t_n} \sum_{i=1}^{N}E_i(t)=\sum_{t=t_0}^{t=t_n} \sum_{i=1}^{N} \lVert \hat{\mathbf{x}}_i(t)-\mathbf{x}_i(t)\rVert^2.
    \label{eq:loss}
\end{equation}

\begin{figure}
    \centering
    \includegraphics[width = 0.4\textwidth]{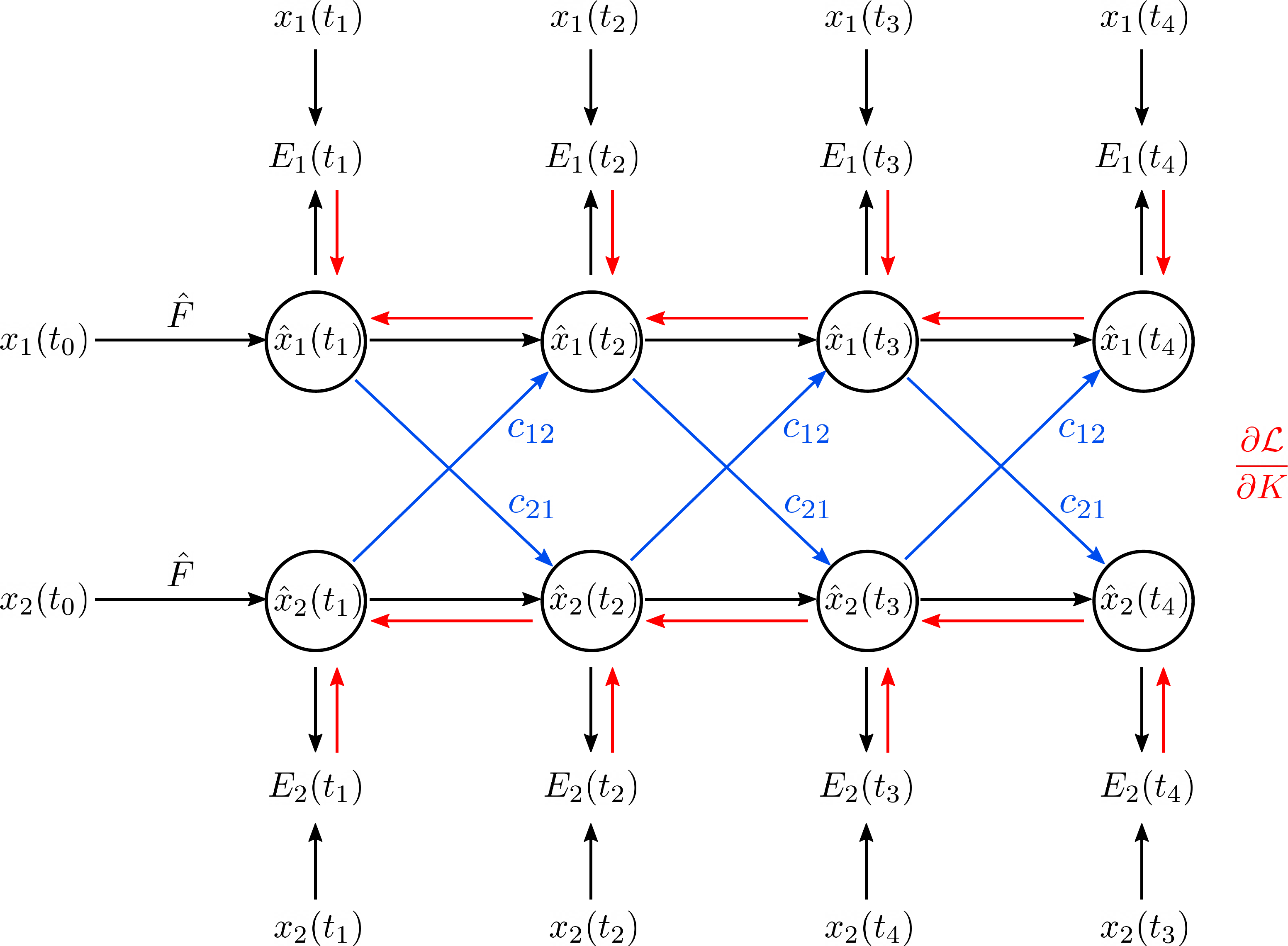}
    \caption{Schematic of the two node modified backpropagation algorithm for calculating weight gradients. Forward pass (blue and black), and backpropagation (red).}
    \label{fig:ModifiedBackprop}
\end{figure}

For regressing network weights, the loss gradient with respect to $\hat{C}$ must be calculated and used to update the estimate of $\hat{C}$,

\begin{equation}
    \medmath{\frac{d \mathcal{L}}{d\hat{C}} = \sum_{t=t_0}^{t=t_n} \sum_{i=1}^{N} \frac{\partial E_i(t)}{\partial \hat{C}} = \sum_{t=t_0}^{t=t_n} \sum_{i=1}^{N} \sum_{d} 2(\hat{\mathbf{x}}^{(d)}_i(t)-\mathbf{x}^{(d)}_i(t)) \frac{\partial \hat{\mathbf{x}}^{(d)}_i(t)}{\partial \hat{C}}.}
    \label{eq:ExpandedLoss}
\end{equation}
This requires the estimation of partial derivatives with respect to the coupling weights $c_{ij}$. The recursive expressions for these derivatives are given in the Appendix. To provide more stability in the regression, the gradient for each step of backpropagation is averaged over the multiple randomly chosen initial inputs $\mathbf{x}(t_0)$. An overview of the modified backpropagation algorithm is given in Figure \ref{fig:ModifiedBackprop}.

\subsubsection*{Decoupling and Retraining}

The use of a mean field approximation to construct an initial local model $\hat{f}$ results in significant \textcolor{edit_col}{errors in the estimation of the vector field for the coupled dimensions.} These errors limit the ability for the algorithm to estimate the true coupling weights as accurate regressions assume that $f =\hat{f}$. By using estimated coupling values $\hat{C}$ to partially decouple calculated values of the vector field, the output space of $\hat{f}$ changes to better reflect the true local dynamics $f$ \textcolor{edit_col}{(i.e. without coupling effects).} To adapt to this change, the model must be retrained on data of the new partially decoupled output space. \textcolor{edit_col}{The process of decoupling the vector field and truncating coupling interaction effects allows for gradual correction of the errors introduced from the mean field approximation of the vector field.}

After attaining an improved estimate of the weighted adjacency matrix $\hat{C}$, the calculated coupling weights can be used to partially decouple the observed node time series according to Equation \ref{eq:decouple}.

\begin{equation}
    \hat{f}(x_i(t_n)) \approx \frac{x_i(t_{n+1})-x_i(t_n)-\delta t \;\sum_{i\neq j}\hat{c}_{ij}g(x_i(t_n),x_j(t_n))}{\delta t}.
    \label{eq:decouple}
\end{equation}

The new calculated value of the vector field $\hat{f}$ can then be used to improve the previously constructed local dynamics model. The final decoupling step is alternated with the backpropagation step to gradually improve both the local model and estimate weighted adjacency matrix $\hat{C}_{ij}$. A much lower model learning rate \textcolor{edit_col}{$\eta' < \eta$} must be used to ensure that the structure of the initial model is preserved. Using a small learning rate also has the added benefit of providing additional stability to the algorithm as it minimises drastic changes in the loss landscape that can affect convergence.

\subsection{Performance Assessment}
We need to introduce several measures to assess performance for the proposed backpropagation method. Firstly, prediction horizons were measured for the regressed dynamical network to assess the overall quality of the combined model inclusive of local dynamics and coupling weights. The prediction horizon for node $i$ was defined as the time $t_{p,i}$ where prediction errors exceeded a set threshold $\epsilon$. \textcolor{edit_col}{The overall prediction horizon $t_p$ for the entire dynamical network is taken as the average prediction horizon across all nodes $i$,}
\begin{equation}
    t_p=\langle t_{p,i}\; |\; \lVert \hat{\mathbf{x}}_i(t_{p,i})-\mathbf{x}_i(t_{p,i}) \rVert^2 <\epsilon \rangle_{i=1}^N.
\end{equation}




\textcolor{edit_col}{The norm of the weight error matrix, $\epsilon_C$ is} used to assess the overall error in the regressed weighted adjacency matrix and is normalised the by the true weight values,

\textcolor{edit_col}{
\begin{equation}
    \epsilon_C = \frac{\lVert \hat{C} -C \rVert ^2}{\lVert C \rVert ^2}.
\end{equation}\\}

\textcolor{edit_col}{
Mutual information $I(X,Y)$ is a nonlinear measure of the correlation between two random variables $X$ and $Y$. In simple terms, describes the amount of information that is shared between any two random variables $X$ and $Y$. Low mutual information implies independence between $X$ and $Y$. Mutual information is usually calculated as a function of marginal and joint probabilities,
\begin{equation}
    I(X,Y) = \int_{\mathcal{X}} \int_{\mathcal{Y}} p_{(X,Y)}(x,y)\ln \left( {\frac{p_{(X,Y)}(x,y)}{p_X(x)p_Y(y)}} \right) \, dydx.
\end{equation}}

\textcolor{edit_col}{
The quality of a local model can be accessed by comparing the mutual information between the true states of a trajectory $\mathbf{x}(t)$, and those predicted by the model $\hat{\mathbf{x}}(t)$ at each point in time $t$. Hence, we can define the average mutual information $I(\tau)$ with respect to some chosen freerun prediction length $\tau$,
\begin{equation}
    I(\tau) = \frac{1}{\tau}\sum_{t=1}^\tau p(\mathbf{x}(t),\hat{\mathbf{x}}(t))\ln \left( {\frac{p(\mathbf{x}(t),\hat{\mathbf{x}}(t))}{p(\mathbf{x}(t))p(\hat{\mathbf{x}}(t))}} \right),
\end{equation}
}
\textcolor{edit_col}{where $p(\mathbf{x}(t)),p(\hat{\mathbf{x}}(t))$ and $p(\mathbf{x}(t),\hat{\mathbf{x}}(t))$ are the joint and marginal probabilities of observing (or predicting for the case of $\mathbf{\hat{x}}$) a given state at time $t$. For each value of $\tau$ in $I(\tau)$, the distributions of $p(\mathbf{x}(t)),p(\hat{\mathbf{x}}(t))$ and $p(\mathbf{x}(t),\hat{\mathbf{x}}(t))$ are calculated empirically using histograms based on observations of $\mathbf{x}(t)$ and $\mathbf{\hat{x}}(t)$ for $t\in[1,\tau]$.}

\textcolor{edit_col}{The function $I(\tau)$ can be interpreted as measure of the prediction performance of the local model for a prediction length of $\tau$ steps. For chaotic systems, $I(\tau)$ of an imperfect model will be a decreasing function. The contour of $I(\tau)$ tracks the collapse of the local model's prediction performance for increasing prediction lengths $\tau$.}

\textcolor{edit_col}{The performance of a trained local model $\hat{f}$ can be assessed by using it to produce freerun predictions of separate system with only a single node (i.e. no external interactions). A score $S$ related to the mutual information between the trajectories from the real local model $\mathbf{x}(t)$ and trained local model $\hat{\mathbf{x}}(t))$ can be calculated using the cumulative area under the average mutual information curve,
\begin{equation}
    S = \int_{1}^{\tau_{max}} I(\tau) d\tau,
\end{equation}
where larger values of $S$ indicate more persistent mutual information and a better predictive local model. }

\textcolor{edit_col}{The mutual information based score $S$ is only calculated with respect to only one component of the predicted time series (i.e. $x$ component). This is sufficient for the systems analysed as all components of the time series ($x,y,z$) feedback to each other. Hence, collapse in the prediction performance of one component will quickly result in poor predictions in all other components. Additionally, $I(\tau)$ for each component cannot be summed together as it assumes that all components evolve independently.}

\textcolor{edit_col}{To assess the mutual information of the combined model of both regressed weights and trained local model, $S$ is calculated as the sum of the scores across all nodes. Whilst nodes are not independent, the sparsity of network connections mean that collapse in the prediction performance of one node does necessarily propagate quickly to all other nodes in the network.}


\section{Results}

Three tests were done to assess the performance of the proposed adapted backpropagation regression algorithm. These were the cases of the small 16 node network, large 64 node network, and networks with non-positive asymmetric coupling weights. \textcolor{edit_col}{All tests were conducted on random networks with each edge existing with probability $p=log(N)/N$ to minimise the likelihood of disjoint networks. Edges are bidirectional in all networks except the asymmetric network case.} \textcolor{edit_col}{Values for the node coupling strengths were normally distributed $(\mu, \sigma)=(0.15,0.02)$. The mean coupling strength $\mu$ is assumed to be non-zero to ensure that regressed values clearly distinguish between the absence and presence of edges.}

Input data of $25000$ steps with an additional $2000$ steps washout were generated with timestep $dt=0.02$ using the standard RK4 algorithm. Washout steps were discarded exclude any transient dynamics. An initial learning rate of $\eta=0.001$ was used for training the feedforward model during initialisation. \textcolor{edit_col}{A lower learning rate of $\eta'=0.0002$ was used for all subsequent retraining iterations with each training period consisting of $30$ epochs. A preset number of $40$ or $80$ refit iterations were chosen before ending regression. However, we note that it may be possible to implement a stopping criterion based on plateauing of either the prediction errors or regressed weight matrix $\hat{C}$.}

\textcolor{edit_col}{Performance measures $(S,t_p)$ for each model were evaluated over 8 randomly chosen initial conditions.To assess the performance of regressed models, each performance measure was compared against a control case where the entire network structure and local model is fully known but initial conditions are perturbed by $\xi \sim N(0,0.005^2)$ where $\xi$ is in scaled units. This method was chosen over the calculation of dynamical invariants such as Lyapunov exponents due to its ease of computation, especially in high dimensional systems such as dynamical networks. Due the chaotic nature of the system's oscillators, small deviations in initial conditions will grow exponentially. Deviations the initial conditions that are sufficiently small are able to track nearby trajectories for some period of time before diverging. Therefore, scores that match those of the control case imply that the regressed models are able to perform almost as well as the case where the system dynamics are fully known, but contaminated with noise in its initial conditions. } 

\subsection{Positive Symmetric Networks}

The fully observable case assumes that all node time series are accessible for regression. Gradient losses $\mathcal{L}$ are calculated at each backpropagation step and used to regress the network adjacency matrix $\hat{C}$. The regression algorithm was run for a total of $30$ refit iterations each consisting of $300$ backpropagation steps followed by decoupling and retraining of the model. A total of 8 randomly generated 16 node networks were tested. The quality of the regressed weights and local model for each network was assessed by running freerun predictions with 16 randomly generated initial conditions.



\begin{figure}[h]
    \centering
    \includegraphics[width = 0.45\textwidth]{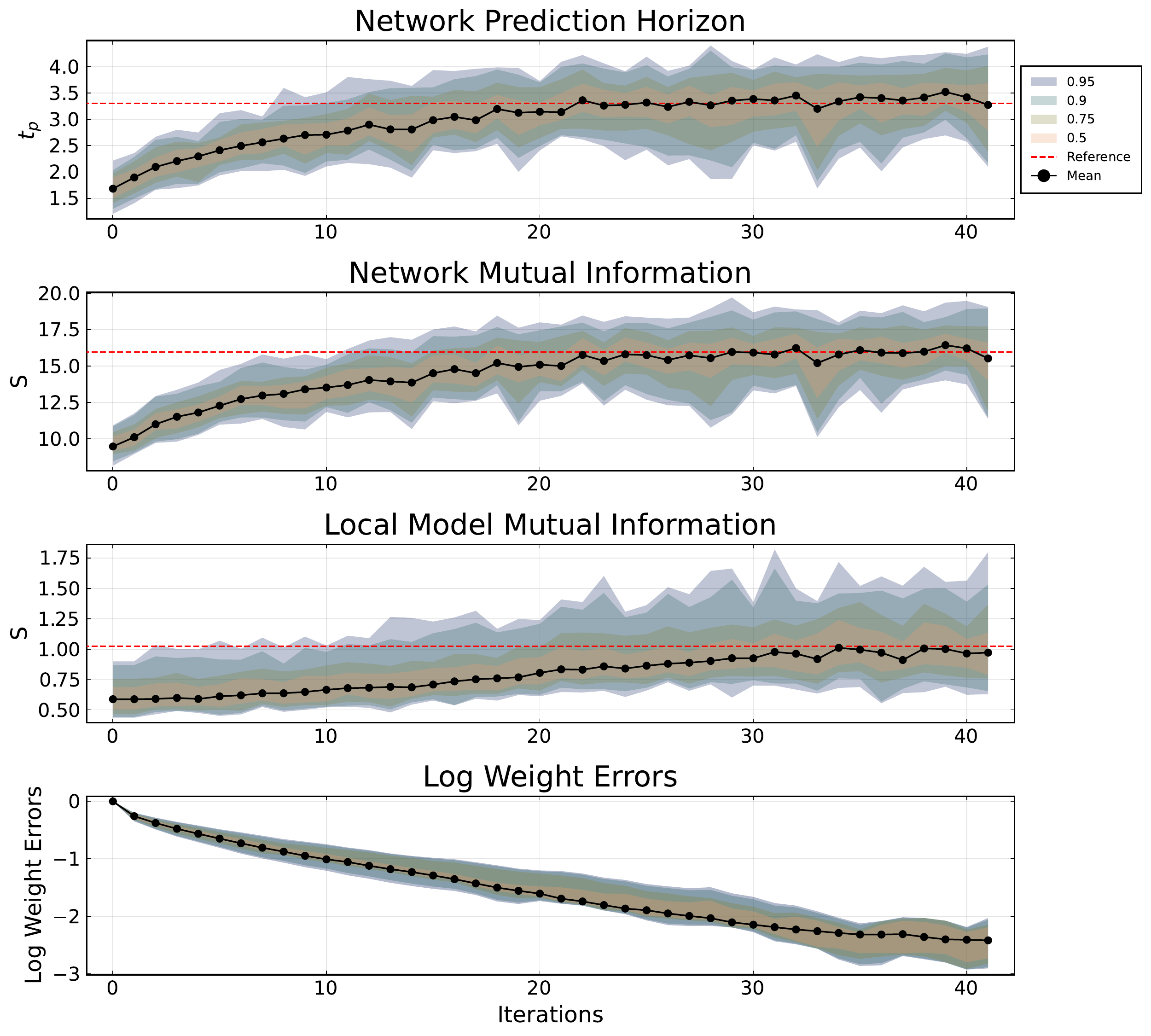}
    \caption{\textcolor{edit_col}{Regression results for the Lorenz 16 node system normally distributed coupling weights $(\mu,\sigma^2)  = (0.15,0.02)$ and connection probability $p$ over 30 iterations. First two data points correspond to the values attained during initialisation.}}
    \label{fig:LorenzCombined}
\end{figure}

\begin{figure}[h]
     \centering
    \includegraphics[width = 0.5\textwidth]{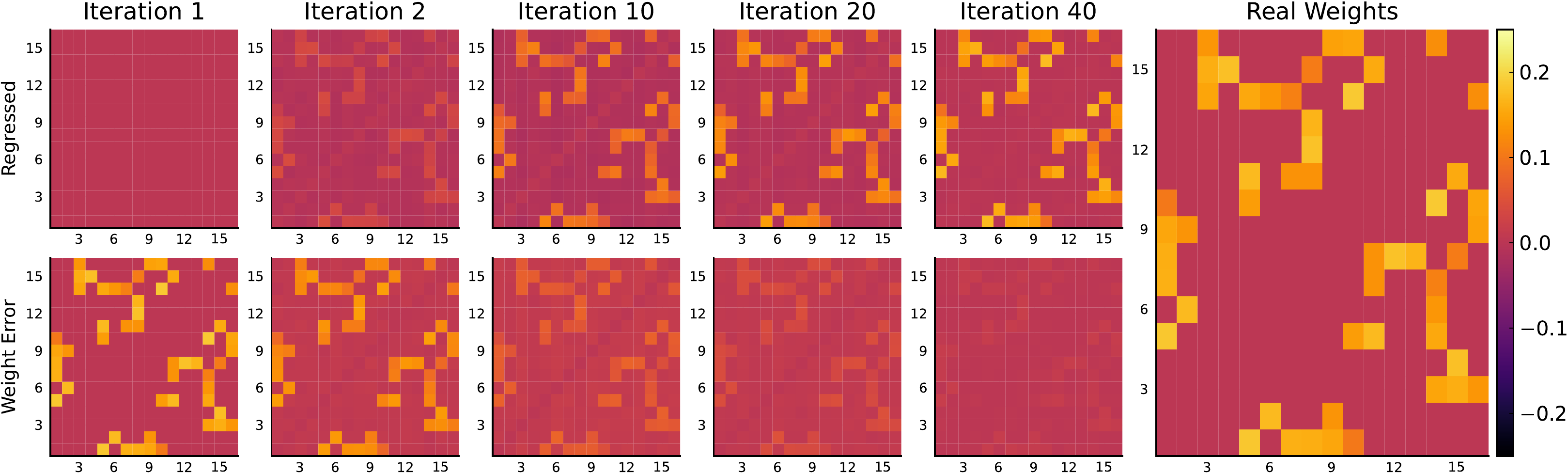}
    \label{fig:LorenzArrays}
     \caption{\textcolor{edit_col}{Regression of Lorenz 16 node network weights (top) at different number of iterations compared to the normalised error in the weights (bottom). True weights (right) are given for comparison showing good agreement with the regressed results.}}
\end{figure}

\begin{figure}
    \centering
    \includegraphics[width = 0.45\textwidth]{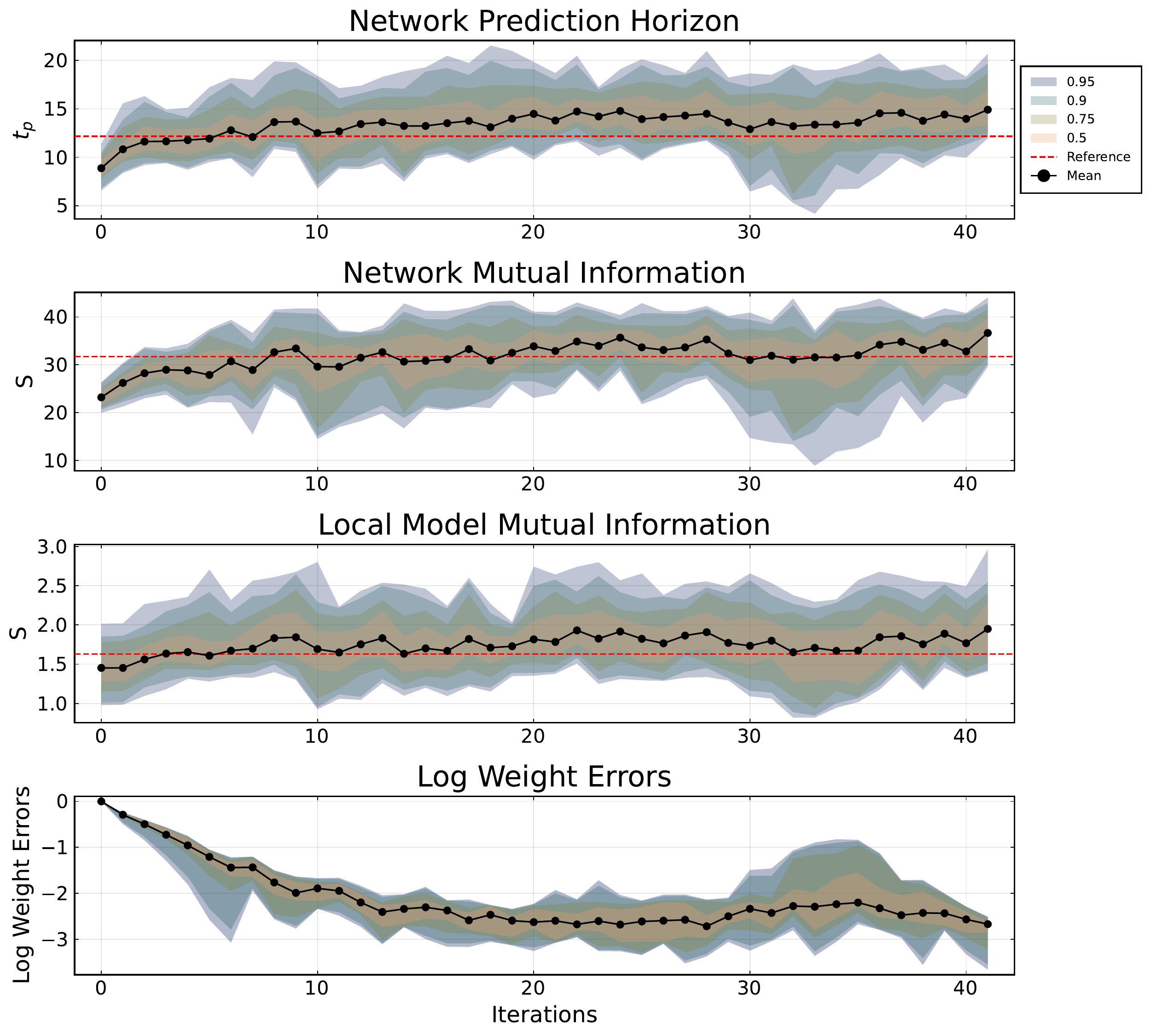}
    \caption{\textcolor{edit_col}{Regression results for the Chua 16 node system operating in the single scroll regime. Coupling weights are normally distributed with $(\mu,\sigma^2)  = (0.15,0.02)$ and connection probability $p$ over 30 iterations. First two data points correspond to the values attained during initialisation.}}
    \label{fig:ChuaCombined}
\end{figure}

\begin{figure}
    \centering
   \includegraphics[width = 0.5\textwidth]{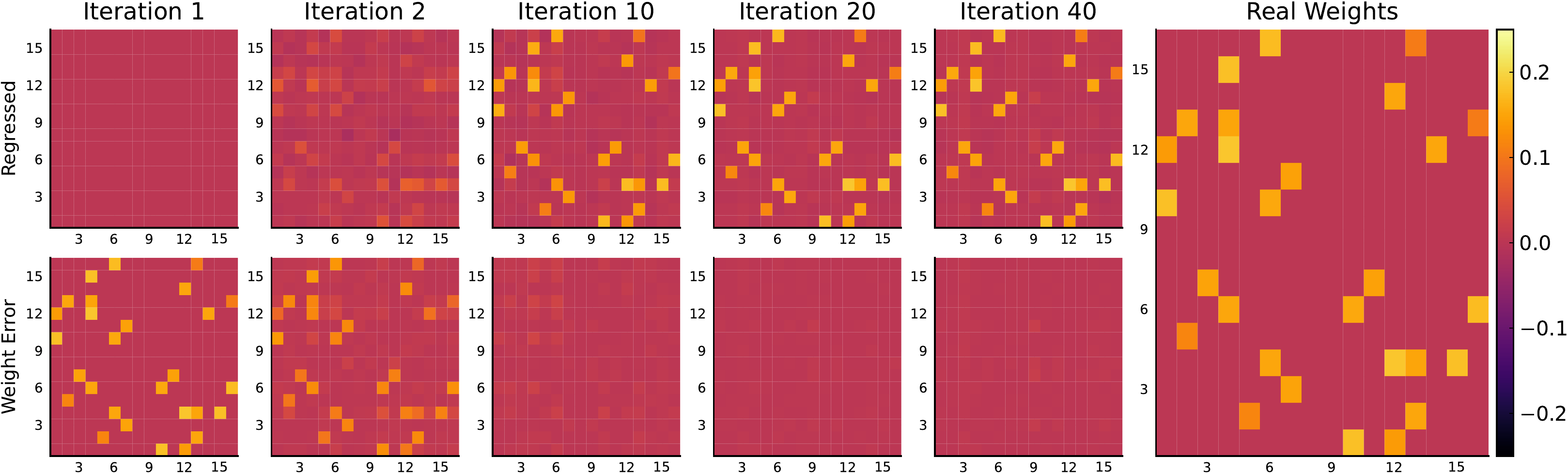}
    \caption{\textcolor{edit_col}{Regression of Chua 16 node network weights (top) at different number of iterations compared to the normalised error in the weights (bottom). True weights (right) are given for comparison showing good agreement with the regressed results.}}
    \label{fig:ChuaArrays}
\end{figure}

\textcolor{edit_col}{The backpropagation method was able to accurately regress the weights of both the Lorenz (Figure \ref{fig:LorenzCombined}) and Chua (Figure \ref{fig:ChuaCombined}) dynamical networks. The average normalised error in the regressed weights for these systems at the end of 30 iterations was 9\% and 7\% for the Lorenz and Chua networks respectively. The constructed trained local model for both systems also showed a large improvement in $S$ with 71\% increase for Lorenz and a 38\% increase for Chua when compared to the initial local model based purely on the mean field approach (i.e. Iteration 0). The mutual information of the combined local model and coupling structure increased by 65\% and 59\%  for the Lorenz and Chua networks respectively. The prediction horizons for both systems were found to increase with each additional iteration. Both Chua and Lorenz networks were found to perform as well or better than the control case where $\xi=0.005$.}

\subsection{Large 64-node Network}

\begin{figure}
    \centering
    \includegraphics[width = 0.45\textwidth]{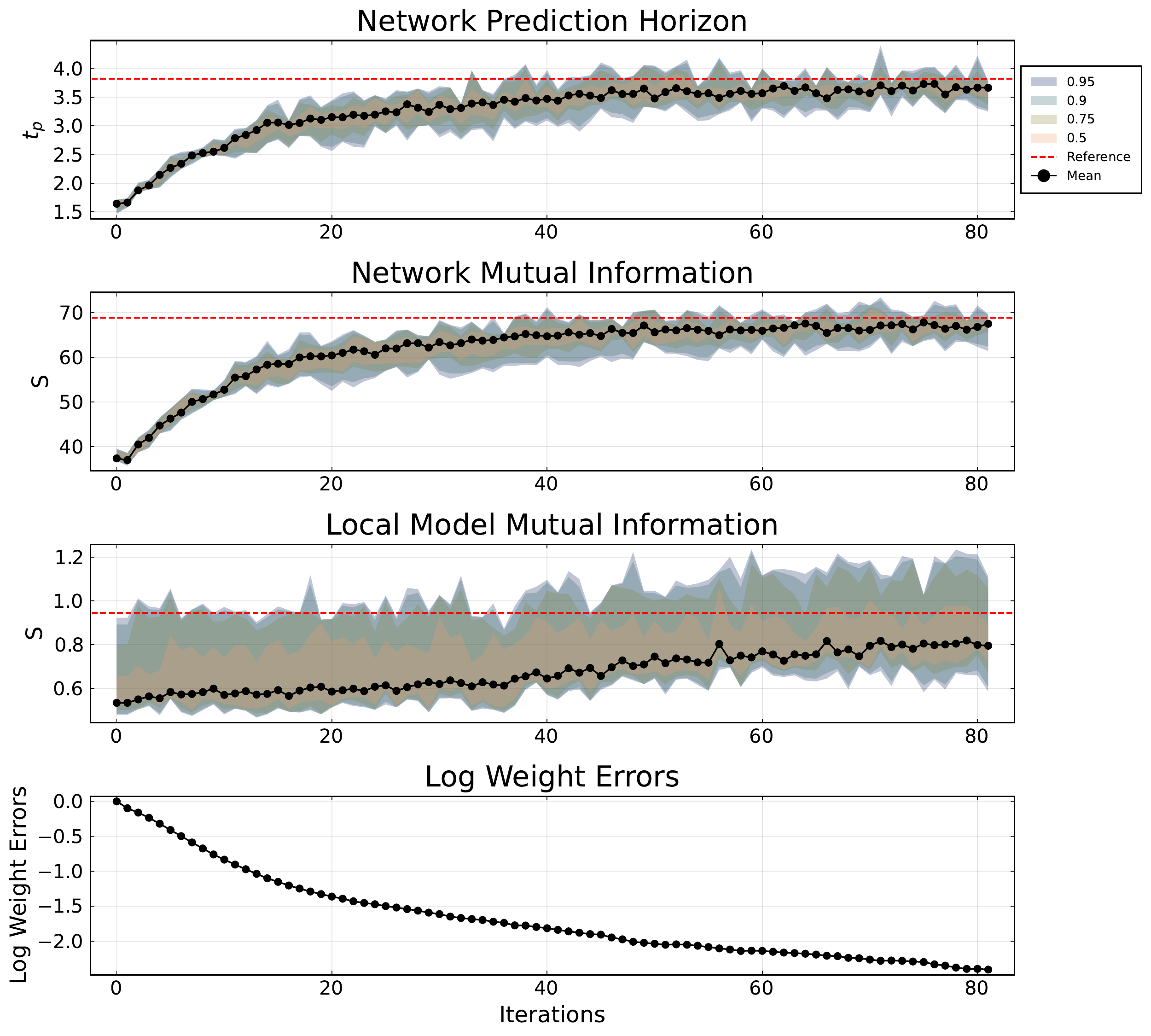}
    \caption{\textcolor{edit_col}{Regression results for the larger Lorenz 64 node system normally distributed coupling weights $(\mu,\sigma^2)  = (0.15,0.02)$ and connection probability $p$ over 80 iterations. First two data points correspond to the values attained during initialisation.}}
    \label{fig:Lorenz64Combined}
\end{figure}

\begin{figure}
    \centering
    \includegraphics[width = 0.5\textwidth]{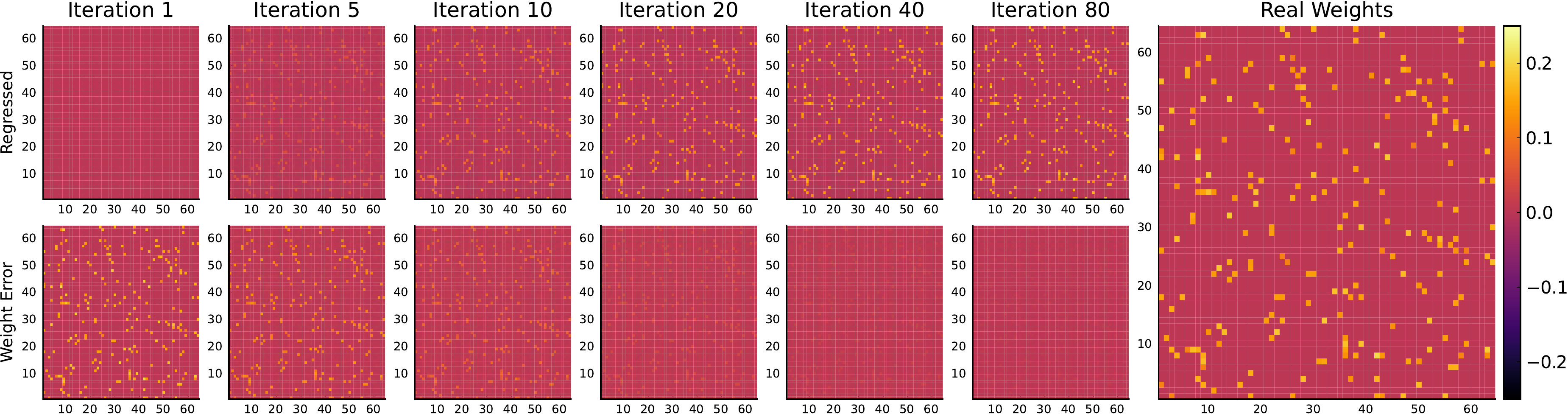}
    \caption{\textcolor{edit_col}{Regression of Lorenz 64 node network weights (top) at different number of iterations compared to the normalised error in the weights (bottom). True weights (right) are given for comparison showing good agreement with the regressed results.}}
    \label{fig:Lorenz64Arrays}
\end{figure}

A single larger system of 64 random coupled Lorenz oscillators was used to test the scalability of the system. 
The algorithm was run for 80 iterations to account for slightly slower convergence. Performance metrics identical to the regular 16 node case was used with the network tested against 16 different initial conditions (see Figures \ref{fig:Lorenz64Combined} and \ref{fig:Lorenz64Arrays}).

We find that the backpropagation method was also able to accurately regress the desired weights of the the coupled system with a final normalised weight error of 9\%. The mutual information of the local model and combined network system was also found to increase by 51\% and 81\% by the end of the iterations respectively.

\subsection{Negative Asymmetric Coupling}

\begin{figure}
    \centering
    \includegraphics[width = 0.45\textwidth]{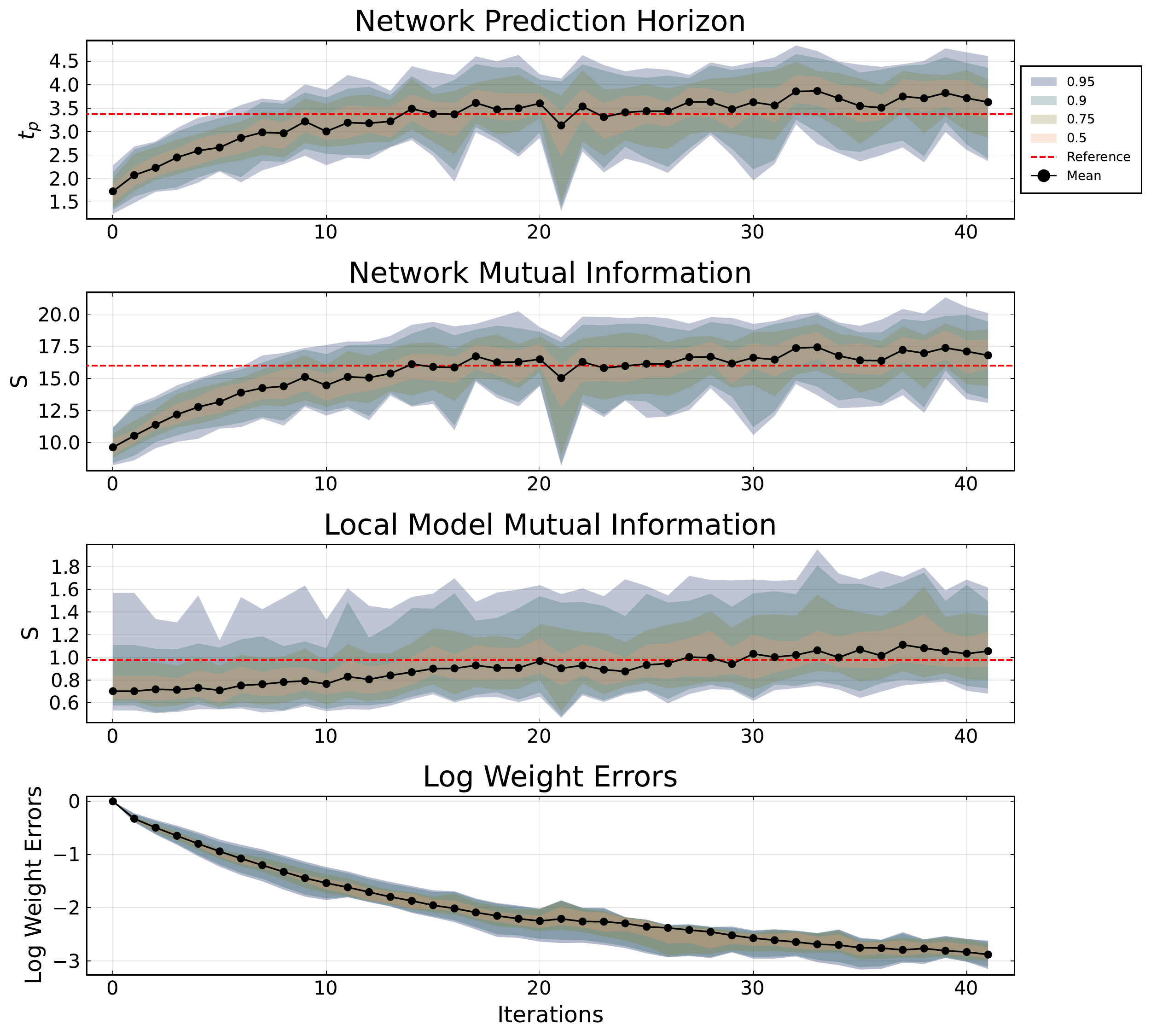}
    \caption{\textcolor{edit_col}{Regression results for the Lorenz 16 node system normally distributed coupling weights $(\mu,\sigma^2)  = (0.15,0.02)$ and connection probability $p$ over 80 iterations. Weights are asymmetric and changed to its negative value with probability of 0.25. First two data points correspond to the values attained during initialisation.}}
    \label{fig:LorenzNegAsymCombined}
\end{figure}

\begin{figure}
    \centering
   \includegraphics[width = 0.5\textwidth]{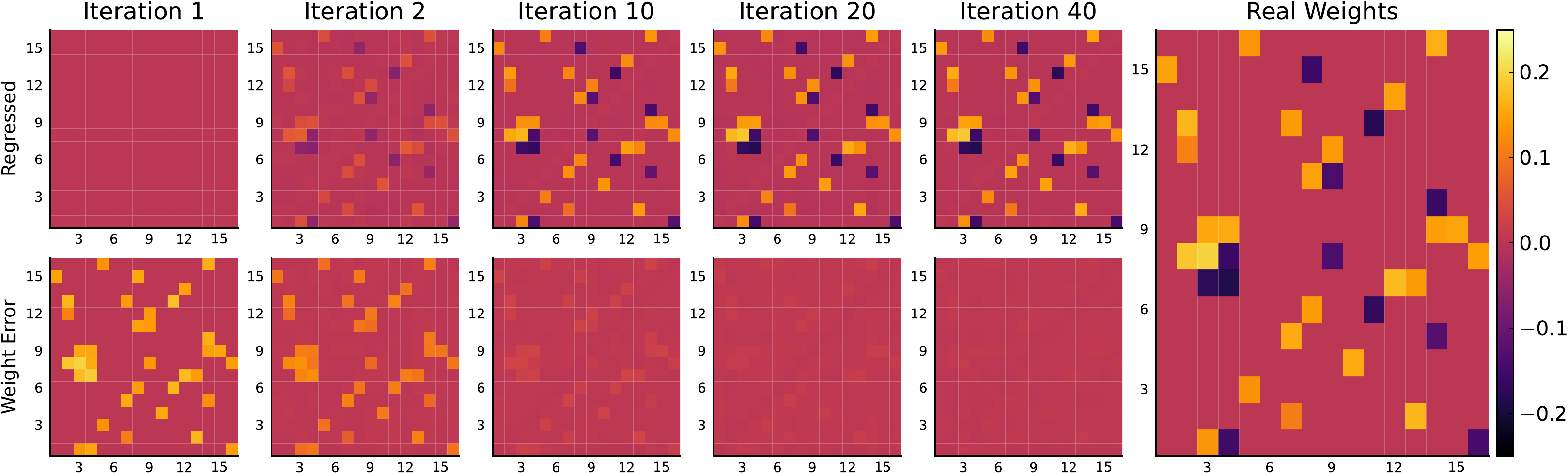}
    \caption{\textcolor{edit_col}{Regression of Lorenz 16 node network with negative and positive weights (top) at different number of iterations compared to the normalised error in the weights (bottom). True weights (right) are given for comparison showing good agreement with the regressed results.}}
    \label{fig:LorenzNegAsymArrays}
\end{figure}

Many dynamical systems assume the case of positive coupling resulting in the attraction of nearby trajectories. In addition to this common case, we consider dynamical systems with both negative and positive coupling where repulsive behaviour between nodes is possible. Additionally, we also consider the effects where the adjacency matrix is not symmetrical. \textcolor{edit_col}{To do this, a 16 node adjacency matrix was generated with each non diagonal entry being non-zero with a probability $p=log(N)/N$ and is not guaranteed to be symmetric. However, each node has a 0.25 probability of being changed to a negative value.} The time series was then generated similarly to the original test case with the Lorenz equations (see Figures \ref{fig:LorenzNegAsymCombined} and \ref{fig:LorenzNegAsymArrays}).

The backpropagation method was found to also work equally well in regressing cases with negative weights and asymmetric coupling. The results for the network prediction horizon, log weight errors and mutual information measures show a similar trend to their symmetric, positively weighted counterparts with a final normalised weight error of 5\%, and an 58\% relative increase in mutual information for the constructed local model. The mutual information of the combined network was found to increase by 75\% when compared to the first iteration.


\textcolor{edit_col}{
\subsection{Noise Effects}
In order to test the robustness of the method, the backpropagation approach was applied to time series from a 16 node Lorenz dynamical network with additive Gaussian noise to simulate measurement noise. The Gaussian noise term $\epsilon_{\xi}(t) \sim N(0,\xi)$ was added after normalising the data to zero mean and unit standard deviation. Because backpropagation regression requires the evaluation of numerical derivatives, the noisy input data was first smoothed using spline regression \cite{green1993nonparametric, reinsch1967smoothing} with a regularisation parameter $\lambda = 10$. The usage of cubic splines also has the advantage or producing relatively smooth derivatives.}

\textcolor{edit_col}{The noise magnitude was varied with $\xi \in [0,0.1]$ with the maximum $\xi$ value corresponding to a signal-noise ratio of 20dB. Backpropagation regression was run for 40 iterations in line with those for the noiseless Lorenz and Chua 16 node systems. Each configuration was repeated for 8 randomly chosen network and initial conditions.}

\begin{figure}
    \centering
    \includegraphics[width = 0.35\textwidth]{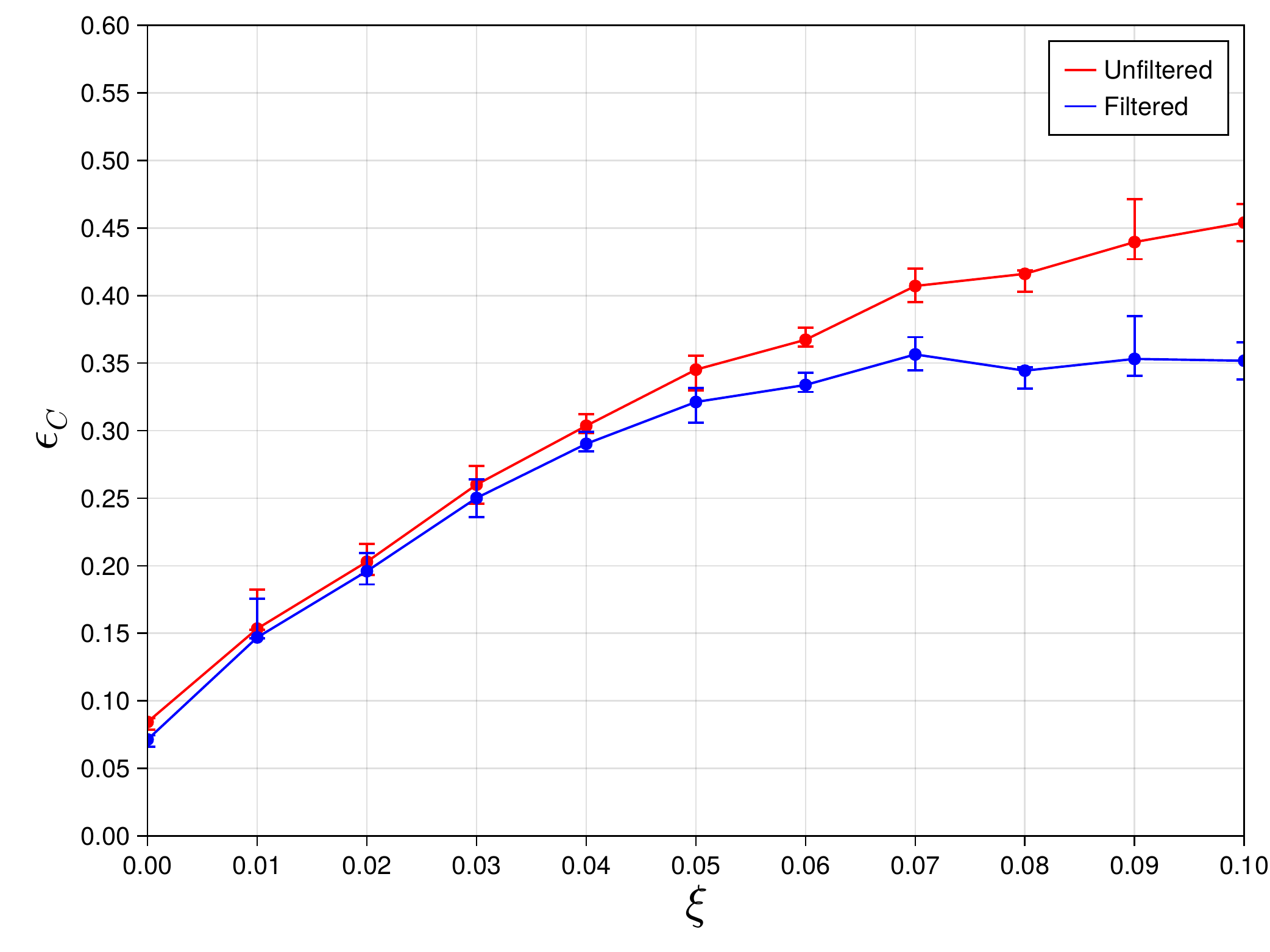}
    \caption{\textcolor{edit_col}{Weight errors $\epsilon_C$ after 40 iterations for backpropagation regression applied to 16 node Lorenz dynamical networks with varying levels of additive Gaussian noise. Each noise configuration was repeated with 8 randomly chosen networks and initial conditions. Error bars correspond to interquartile range of weight errors. Weights with a magnitude less than a threshold 0.04 are truncated and set to zero. The threshold value selected selected from the histogram distribution of all weights.}}
    \label{fig:LorenzNoisyWeights}
\end{figure}

\textcolor{edit_col}{The effect of noise hampers performance and results in the addition of spurious incorrect weights during regression (see Figure \ref{fig:LorenzNoisyWeights}). To address this, a threshold value of 0.04 is selected to filter the weights following regression. Regressed weights that do not exceed this value are set to 0. The threshold value is visually selected from the histogram distribution of weights. As the network weights are normally distributed with non-zero mean, $c_{ij} \sim N(0.15,0.02)$, true weights will tend to have regressed value much greater than zero. In contrast, noise effects will yield to low valued spurious weights. This results in a bimodal distribution of weights. Hence, the threshold value can be selected as a value that distinguishes between the two modes.}

\textcolor{edit_col}{The performance of the regression algorithm weights was found to be relatively robust against noise. For low noise levels $\xi \leq 0.04$, the effects of noise appear follow a linear relationship (see Figure \ref{fig:LorenzNoisyWeights}). The increase in weight errors $\epsilon_C$ slow for higher levels of noise. In the latter case, the weight filtration was found to be very effective in removing spurious edge weights allowing a further decrease in weight errors. However, we found that higher values of noise with $\xi \geq 0.15$ resulted in instability of the regression algorithm and divergence in the weight errors.}

\begin{figure}
    \centering
    \includegraphics[width = 0.35\textwidth]{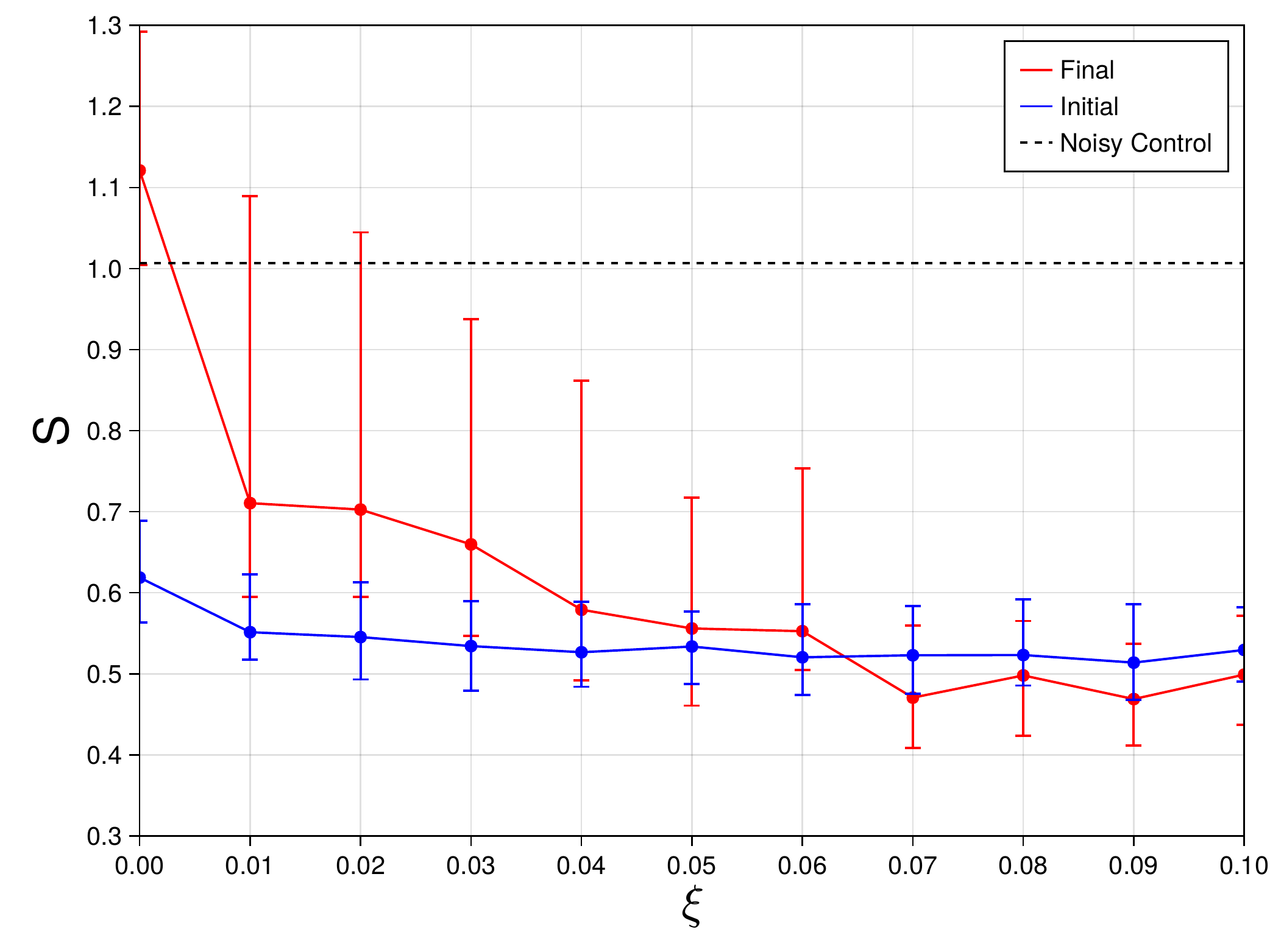}
    \caption{\textcolor{edit_col}{Mutual information of the local models for configurations with increasing observational noise $\xi$. Model scores are compared at the beginning (blue) and end (red) of the 40 refit iterations. Mutual information is compared against the control case where the exact model known with perturbed initial conditions.}}
    \label{fig:LorenzNoisyMI}
\end{figure}

\textcolor{edit_col}{The effect of noise was found to significantly affect the quality of the local dynamics models (see Figure \ref{fig:LorenzNoisyMI}). For small noise levels $\xi \leq 0.01$, the approximations $\hat{f}$ of the local dynamics was relatively stable. However, there was a large difference in prediction performance between the noiseless and noisy case. Noise effects are more detrimental for $\xi > 0.01$ where the long run predictions of the resulting models no longer conform to the attractor of the original Lorenz system. In several instances, it was found that the models slowly drifted and converged to a fixed point. For $\xi \geq 0.07$, the observational noise results in a local model that is worse than the initial mean field estimate. }

\textcolor{edit_col}{The poor learning of the local model can be attributed to the combined effects of noise changing the temporal structure of the time series and the relatively low sampling rates of the time series. An integration step of $dt = 0.02$ corresponds to approximately 35 points per oscillation. Each oscillation conists of a short transient extrema. Additive noise result in over or underestimation of these peaks, which can have a detrimental effect on calculated vector field. The noise tests were repeated at a lower integration step of $dt = 0.002$, which was smoothed and downsampled to $dt = 0.02$. This configuration allowed a better treatment of the additive noise and yielded more robust results.}

\textcolor{edit_col}{Overall, we find that the performance of this method under noise is comparable to other methods used to analysed dynamical networks. For local dynamics, the attractor of the regressed local model dynamics began to collapse at $\xi = 0.05$ corresponding to a signal-noise ratio of 26dB. This is comparable to the test results conducted by \cite{gao2022autonomous} where the methods SINDy, ARNI and Gao's method exhibit large decreases in performance at approximately 55dB, 40dB and 25dB respectively.}

\textcolor{edit_col}{
\section{Method Limitations}
There are several data and system requirements needed for the application of the backpropagation regression algorithm. Firstly, the backpropagation regression method requires networks that are not fully synchronised. This is because differences in the trajectories of coupled nodes are used to decouple observed node signals to infer links. Networks operating in the fully synchronised regime will have fully coincident trajectories and cannot be used to infer weights.
}

\textcolor{edit_col}{Secondly, reliance of a mean field approximation during the initialisation stage assumes that node trajectories do not deviate too much from the attractor of the local dynamics of the uncoupled single node case. Small deviations will inevitably occur, and is accounted for during the decoupling and retraining stage of each iteration. However, deviations that are too large may cause convergence problems.}

\textcolor{edit_col}{Thirdly, there must be sufficient full node state observations with trajectories that fully explore the state space for the regression to converge and produce a good local model. The method was found to be convergent with time series lengths as low as 5000 time steps with $dt = 0.02$, excluding an initial washout period of $2000$ time steps to account for dynamical transients.}

\section{Conclusion}
Dynamical networks are an interesting and useful framework to analyse complex systems with multiple interacting components. However, application of this framework is usually difficult due to lack of information pertaining to features such as the connectivity structure and local dynamics. As a result, there is a need for a method to extract these features purely from observed node states.

In this paper, we propose an adapted backpropagation method for inferring the local model and connectivity structure of dynamical networks purely from node time series. This method draws inspiration from the backpropagation through time (BPTT) algorithm commonly used to train RNNs and adapts it for the purpose of regressing couplings weights in a dynamical network. Regressed values for the network coupling are then used to decouple the observed node time series and improve the construction of the local model without coupling effects. The two steps of backpropagation weight regression and decoupling alternate until sufficient convergence is achieved yielding a final local model and regressed coupling weights.

\textcolor{edit_col}{The method was tested on dynamical networks of nodes of identical chaotic oscillators with successful results. Two types of chaotic oscillators, Lorenz and single-scrolled Chua, were tested. The backpropagation regression method was able to accurately identify the coupling weights down to a margin of 7\%. The regressed models were also compared to a control case where the true local model is known, but whose initial states are perturbed by a small error. It was found that the trained local model was comparable to the control case. Combined, the regressed weights and constructed local model were able to achieve respectable prediction horizons for the 16 node dynamical network compared to a control case. Despite its poor scalability, the method was also able to accurately reconstruct the local model and coupling weights of a larger 64 node network with similar results to the smaller 16 node counterparts. Additional testing was also done on dynamical networks with negative and asymmetric coupling with similar performance. }

\textcolor{edit_col}{The robustness of the method was also tested against noise and network heterogeneity (see Appendix). The algorithm was found to be convergent for small to moderate amounts of observational noise. However, larger noise levels resulted in poor local models and potential instability. In contrast, the presence of network heterogeneity did not significantly affect the performance of the backpropagation regression algorithm.}

\textcolor{edit_col}{We acknowledge that the method proposed in this paper is a proof of concept that, whilst functional, has multiple areas of further research and potential refinement. Firstly, there remains further work on the refinement of the speed and complexity of the algorithm to be more scalable.
Secondly, there remains a need for further testing the on the effect of hyperparameters and initialisation on algorithm performance. Thirdly, it would be beneficial to investigate the performance of this method on different dynamical network structures. There also remains the possibility of investigating how other machine learning methods used in recurrent neural networks may be applied to the backpropagation regression algorithm.}

\section*{Acknowledgments}
E.T. is supported by a Robert and Maude Gledden Postgraduate Research Scholarship and Australian Government RTP Scholarship at The University of Western Australia. M.S. and D.C.C. acknowledge the support of the Australian Research Council through the Centre for Transforming Maintenance through Data Science (grant number IC180100030), funded by the Australian Government.

\ifCLASSOPTIONcaptionsoff
  \newpage
\fi

\bibliographystyle{IEEEtran}
\bibliography{references}

%

\begin{IEEEbiographynophoto}{Eugene Tan}
is a Robert and Maude Gledden scholar and PhD candidate at the University of Western Australia. He completed his B.Sc and MPE in mechanical engineering at the University of Western Australia.
\end{IEEEbiographynophoto}
\begin{IEEEbiographynophoto}{D\`{e}bora Corr\^{e}a}
completed her PhD degree in Applied Physics with a focus on Computational Physics at the University of Sao Paulo (USP). She is currently a lecturer at the University of Western Australia
\end{IEEEbiographynophoto}
\begin{IEEEbiographynophoto}{Thomas Stemler}
completed his Ph.D. at the Technische Universit\"{a}t Darmstadt. He is currently a senior lecturer at the University of Western Australia.
\end{IEEEbiographynophoto}
\begin{IEEEbiographynophoto}{Michael Small}
is a professor and CSIRO-UWA Chair in Complex Engineering Systems at the University of Western Australia. He completed his Ph.D. in applied mathematics focusing on nonlinear dynamical systems and time series analysis at The University of Western Australia.
\end{IEEEbiographynophoto}






\end{document}